\documentclass{article}
\usepackage{amssymb,amsmath,bm,geometry,graphics,graphicx,color,hyperref}
\usepackage{color}
\newtheorem{theo}{Theorem}

\newtheorem{cor}{Corollary}
\newtheorem{rem}{Remark}

\def\proof{\noindent\underline{Proof}\quad}
\def\QED{\mbox{~$\Box{~}$}}
\def\no{\noindent}

\def\RR{\mathbb{R}}

\def\dd{\mathrm{d}}
\def\mat{\begin{pmatrix}}
\def\endmat{\end{pmatrix}}
\def\bfb{\bm{b}}
\def\bfc{\bm{c}}
\def\bfe{\bm{e}}

\def\bfeta{\bm{\eta}}

\def\bfzero{\bm{0}}
\def\bfgamma{\bm{\gamma}}

\def\bfq{\bm{q}}
\def\bfp{\bm{p}}
\def\bfw{\bm{\omega}}
\def\P{{\cal P}}
\def\I{{\cal I}}
\def\H{{\cal H}}
\def\M{{\cal M}}

\def\L{{\cal L}}
\def\sech{\mathrm{sech}}

\title{Spectrally Accurate Energy-preserving Methods for the Numerical Solution of the ``Good" Boussinesq Equation}
\author{Luigi Brugnano$^*$,~ Gianmarco Gurioli$^*$,~ Chengjian Zhang$^\dag$ \\
~\\
\small $^*$ Dipartimento di Matematica e Informatica ``U.\,Dini'', Universit\`a di Firenze, 50134 Firenze, Italy\\
\small $^\dag$ School of Mathematics and Statistics, Huazong University of Science and Technology\\
\small  Wuhan 430074, Hubei, China
}

\begin{document}
\maketitle

\begin{abstract} In this paper we study the {\em geometric} solution of the so called ``good'' Boussinesq equation.
This goal is achieved by using a convenient space semi-discretization, able to preserve the corresponding Hamiltonian
structure, then using energy-conserving Runge-Kutta methods in the HBVM class for the time integration. 
Numerical tests are reported, confirming the effectiveness of the proposed method. 

\medskip
\no {\bf Keywords:} ``good'' Boussinesq equation, Hamiltonian PDEs, energy-conserving methods, Hamiltonian Boundary Value Methods, HBVMs, spectral methods, blended iteration.

\medskip
\no {\bf MSC:} 65P10, 65L05, 65M70.
\end{abstract}

\section{Introduction}
We here consider the efficient numerical solution of the ``good" Boussinesq equation,
\begin{equation}\label{GB1}
w_{tt}(x,t)=-w_{xxxx}(x,t)+w_{xx}(x,t)+(w^2(x,t))_{xx}, \qquad (x,t)\in [a,b]\times[0,\infty),
\end{equation}
commonly used to describe small amplitude long waves propagation on the surface of shallow water. It is for this reason that the equation is often considered in several physical contests, such as ocean and coastal engineering (as stressed, e.g., in \cite{Abd,Moh}). Moreover, the equation provides a balance between dispersion and nonlinearity that may lead either to the existence of solitons, or blowup solutions \cite{Manorajan,Manoranjan2,Zhang,Nguyen,Xiao,Runz,Sachs}. 

For sake of brevity, when not required by the context, we shall skip the arguments $(x,t)$ for the wave function $w$ and its derivatives.
The equation \eqref{GB1} is completed with initial conditions
\begin{equation}\label{inico1}
w(x,0)=w_0(x),\quad w_t(x,0)=g(x),\qquad x\in[a,b],
\end{equation}
and periodic boundary conditions. 
We shall also assume that the wave has the linear invariant
\begin{equation}\label{w0}
\L[w](t) := \int_a^b w(x,t)\dd x \equiv const = \int_a^b w_0(x)\dd x, \qquad \forall t\ge0,
\end{equation}
which, in turn, upon regularity assumptions on the initial data, implies that
\begin{equation}\label{g0}
\int_a^bg(x)\dd x=0.
\end{equation}
In fact, since $\L[w]$ is conserved, one has:\footnote{As usual, the $\dot{~}$ will denote the derivative w.r.t. $t$.}
$$\dot\L[w](0) =\int_a^b w_t(x,0)\dd x =\int_a^bg(x)=0.$$
Hereafter, we shall assume that $w_0(x)$ and $g(x)$ are such that the solution is regular enough, as a periodic function on $[a,b]$, for all $t\ge0$. In order to simplify the
arguments, it is customary to define a shifted variable, 
\begin{equation}\label{u2w}
u(x,t) = w(x,t) + \frac{1}{2}, 
\end{equation}
which transforms equation \eqref{GB1} into the simpler form:\footnote{Hereafter, for sake of brevity we shall omit the arguments of the functions, unless they are needed.}
\begin{equation}\label{GB2}
u_{tt}=-u_{xxxx}+(u^2)_{xx}, \qquad (x,t)\in [a,b]\times[0,\infty),
\end{equation}
with initial conditions (see (\ref{inico1}) and (\ref{u2w}))
\begin{equation}\label{inico2}
u(x,0)=w_0(x)+\frac{1}{2} =: u_0(x),\qquad u_t(x,0)=w_t(x,0)=g(x),\qquad x\in[a,b],
\end{equation}
and periodic boundary conditions. Moreover, since $u$ and $w$ differs by a constant, because of (\ref{w0}) and (\ref{g0}), one has: %we also assume:
\begin{equation}\label{u0}
\L[u](t) = \int_a^b u(x,t)\dd x \equiv \int_a^b u_0(x)\dd x, \qquad \forall t\ge0.
\end{equation}

The numerical solution of \eqref{GB1} or \eqref{GB2} has been developed along different directions, ranging from the pseudo-spectral or splitting approach \cite{Frutos,Chen,Zhang2,Zhang3,Cheng,Cai,Uddin}, up to  finite-difference and finite-element schemes \cite{bra,El,Ortega,Ismail1,pani}, as well as structure-preserving methods \cite{Chen2,huang,Zeng} and energy-preserving methods \cite{Jiang}. In particular, \cite{Yan} and \cite{Yan2} consider an energy-conserving strategy based on the HBVMs for the ``good'' Boussinesq and the improved Boussinesq equation, respectively, while a second order symplectic method preserving the energy and the momentum is considered in \cite{Aydin}.

Hereafter we shall focus on the {\em geometric} numerical solution of the simpler form (\ref{GB2}),
where, by {\em geometric} it is meant that we will provide a numerical solution able to retain important geometric
properties of the continuous one. In particular, we shall see that the equation (\ref{GB2}) has a Hamiltonian structure, which can
be preserved by a suitable space semi-discretization. The time integration will be then performed by using energy-conserving
methods in the HBVMs class \cite{LIMbook2016,BI2018,BIT2009,BIT2010,BIT2012,BIT2015,IaPa2007,IaPa2008,IaTr2009}, and this will allow us to retain many geometric properties of the solution, as later specified: as matter of fact, this paper follows a systematic study of the application of HBVMs for efficiently solving Hamiltonian PDEs \cite{BBFCI2018,BFCI2015,BGS2018,LIMbook2016,BI2018,BZL2018}. In particular, we shall derive a very efficient solution procedure, which fully exploit the particular structure of the problem, in order to define a spectrally accurate numerical method, both in space and time, able to retain relevant geometric properties.

With this premises, the structure of the paper is as follows: in Section~\ref{HF} we study the Hamiltonian formulation of (\ref{GB2}); in Section~\ref{SD}
we study a convenient space semi-discretization; in Section~\ref{HBVM} we sketch the main facts concerning HBVMs, along with their efficient
implementation; in Section~\ref{NT} we report some numerical tests aimed at assessing the geometric properties of the resulting method; 
at last, a few concluding remarks will be given in Section~\ref{end}.
 
\section{Hamiltonian formulation}\label{HF}
The equation (\ref{GB2}) can be recast as an Hamiltonian system of two partial differential equations, as follows:
\begin{equation}\label{uvform}
u_t = v_x, \qquad v_t = -u_{xxx} + (u^2)_x, \qquad (x,t)\in[a,b]\times[0,\infty),
\end{equation}
with initial conditions (see (\ref{inico2}))
\begin{equation}\label{uv0}
u(x,0)= u_0(x),\qquad v(x,0) = \int_a^x g(s)\dd s =: v_0(x),\qquad x\in[a,b],
\end{equation}
and periodic boundary conditions. In fact, $u$ has periodic boundary conditions and, because of (\ref{u0}) and (\ref{uvform}), one also has:\footnote{This technical detail is often overlooked in the literature.}
$$v(b,t)-v(a,t) = \int_a^b v_x(x,t)\dd x =\int_a^b u_t(x,t)\dd x = \dot\L[u](t) = 0, \qquad \forall t\ge0.$$
In particular, one obtains that the system (\ref{uvform}) can be formally written in a more compact way as:
\begin{equation}\label{uvform1}
\begin{pmatrix}  u_t \\ v_t \end{pmatrix}=|J_2|\otimes \partial_x ~ \delta \H [u,v],\end{equation}
with 
\begin{equation}\label{J2dH}
J_2=\begin{pmatrix} &1\\-1\end{pmatrix} \qquad \mbox{and}\qquad
\delta\H[u,v]= \begin{pmatrix}  \delta_u \H[u,v] \\ \delta_v\H[u,v] \end{pmatrix}
\end{equation}
the vector of the functional derivatives of the Hamiltonian functional
\begin{equation}\label{Hfun}
\H[u,v]=\frac{1}{2}\int_a^b\left(v^2+\frac{2}{3}u^3+u_x^2 \right)\dd x =: \int_a^b L(v,u,u_x)\dd x.
\end{equation}
In fact, one has:
\begin{eqnarray*}
\delta_u\H [u,v] &=& \left( \partial_u- \partial_x\partial_{u_x}\right) L(v,u,u_x) ~=~u^2-u_{xx},\\
\delta_v\H [u,v]&=&\partial_vL(v,u,u_x)~=~v.
\end{eqnarray*}
Therefore, the ``good" Boussinesq equation is an instance of a second order {\em Hamiltonian PDE}.
 
\begin{theo}\label{Hcost} 
Assume that the solution of (\ref{uvform}) is $C^3[a,b]$ as a periodic function. Then, the Hamiltonian functional (\ref{Hfun}) is constant along the solution of (\ref{uvform}).
\end{theo}
\proof In fact, by the hypotheses on $u$ and $v$, one has:
\begin{eqnarray*}
\dot\H[u,v]&=&\frac{1}{2}\int_a^b\left(2vv_t+2u^2u_t+2u_xu_{xt}\right)\dd x=\int_a^b\left[v(-u_{xxx}+2uu_x)+u^2v_x+u_xv_{xx}\right]\dd x\\
&=&\int_a^b\left(2uvu_x+u^2v_x+u_xv_{xx}-vu_{xxx}\right)\dd x=\int_a^b \left(u^2v+v_xu_x-vu_{xx}\right)_x\dd x\\
&=&\left[u^2v+v_xu_x-vu_{xx}\right]_{x=a}^{x=b}=0,
\end{eqnarray*}
because of the periodicity in space of the functions $u$ and $v$, as well as their derivatives w.r.t. $x$.\QED

\medskip
In addition, we can consider the following quadratic functional
\begin{equation}\label{M}
\M [u,v]=\int_a^b uv\,dx,
\end{equation}
corresponding to the {\em momentum} (or the {\em impulse}), for which the following result holds true.

\begin{theo}\label{Mcost}
In the same hypotheses of Theorem~\ref{Hcost}, the quadratic functional (\ref{M}) is constant along the solution of (\ref{uvform}).
\end{theo}
\proof In fact, using arguments similar to those used in the previous theorem, one has:
\begin{eqnarray*}
\dot\M[u,v]&=&\int_a^b\left(u_tv+uv_t\right)\dd x=\int_a^b\left[vv_x+u(-u_{xxx}+2uu_x)\right]\dd x\\
&=&\int_a^b\left(\frac{1}{2}v^2+\frac{2}{3}u^3+\frac{1}{2}u_x^2-uu_{xx}\right)_x\dd x
=\left[\frac{1}{2}v^2 +\frac{2}{3}u^3+\frac{1}{2}u_x^2-uu_{xx}\right]_{x=a}^{x=b}=0,
\end{eqnarray*}
by virtue of the periodicity in space of the functions $u$ and $v$, as well as their derivatives w.r.t. $x$.\QED

\medskip
It is worth mentioning that, besides (\ref{u0}), also the functional (see (\ref{uv0}))
\begin{equation}\label{v0}
\L[v](t) := \int_a^b v(x,t)\dd x, \qquad \forall t\ge0,
\end{equation}
is conserved. In fact, one has:
\begin{equation}\label{v00}
\L[v](0)=  \int_a^b v_0(x) \dd x
\end{equation}
and
$$\dot\L[v] = \int_a^b v_t\dd x = \int_a^b \left(-u_{xxx}+(u^2)_x\right)\dd x= \int_a^b \left(-u_{xx}+u^2\right)_x\dd x = \left[ -u_{xx} + u^2\right]_{x=a}^{x=b} = 0,$$
because of the periodicity of $u$ (and its space derivatives).

\medskip
In particular, conserving $\L[u]$, $\H[u,v]$, $\M[u,v]$, $\L[v]$ in (\ref{u0}), (\ref{Hfun}), (\ref{M}), and (\ref{v0}), represents the relevant {\em geometric properties} of the solution we are interested in, which we shall try to reproduce in the discrete approximation.

\section{Space discretization}\label{SD}
We now discretize the space variable along the following orthonormal basis for periodic $L^2[a,b]$ functions:
\begin{equation}\label{cs}
c_j(x)=\sqrt{\frac{2-\delta_{j0}}{b-a}}\cos\left(2\pi j\frac{x-a}{b-a}\right),\quad j\ge 0,\qquad
s_j(x)=\sqrt{\frac{2}{b-a}}\sin\left(2\pi j\frac{x-a}{b-a}\right),\quad ~~j\ge1.
\end{equation}
In fact, one verifies that, for all allowed indexes $i,j$: 
\begin{equation}\label{orto}
\int_a^b c_i(x)c_j(x)\,dx=\delta_{ij}=\int_a^b s_i(x)s_j(x)\,dx,\qquad \int_a^b s_i(x)c_j(x)\,dx=0.
\end{equation}
Consequently, for suitable time dependent coefficients $\alpha_j(t), \beta_j(t), \xi_j(t),\eta_j(t)$, the following expansions are derived:
\begin{eqnarray}\nonumber
u(x,t)&=&\alpha_0(t)c_0(x)+\sum_{j\ge 1}\left( \alpha_j(t)c_j(x)+\beta_j(t)s_j(x)\right),\\[-2mm] \label{temp}\\[-2mm] \nonumber
v(x,t)&=&\xi_0(t)c_0(x)+\sum_{j\ge 1}\left( \xi_j(t)c_j(x)+\eta_j(t)s_j(x)\right).
\end{eqnarray}
One easily verifies the following result.

\begin{theo}\label{u0v0}
In order to conserve $\L[u]$ and $\L[v]$, see  (\ref{u0}) and (\ref{v0})-(\ref{v00}), respectively, in the expansions (\ref{temp}) one must have:
$$%\begin{equation}\label{hu0}
\alpha_0(t)c_0(x) \equiv \frac{1}{b-a} \int_a^b u_0(x)\dd x =: \hat u_0, \qquad \xi_0(t)c_0(x) \equiv \frac{1}{b-a} \int_a^b v_0(x)\dd x =: \hat v_0.
$$%\end{equation}
\end{theo}
\proof The statements easily follows from the fact that $\int_a^b c_j(x)\dd x =\int_a^bs_j(x)\dd x=0$, for all $j=1,2,\dots.\,\QED$

\bigskip
As a consequence, the previous expansions (\ref{temp}) becomes:
\begin{eqnarray}\nonumber
u(x,t)&=&\hat u_0+\sum_{j\ge 1}\left( \alpha_j(t)c_j(x)+\beta_j(t)s_j(x)\right) \equiv \hat u_0+\bfw(x)^\top\bfq(t),\\[-2mm] \label{uv}
\\[-2mm] \nonumber \label{v}
v(x,t)&=&\hat v_0+\sum_{j\ge 1}\left( \xi_j(t)c_j(x)+\eta_j(t)s_j(x)\right) \equiv \hat v_0+\bfw(x)^\top\bfp(t),
\end{eqnarray}
having set the infinite vectors
\begin{equation}\label{wqp}
\bfw(x)=\begin{pmatrix} s_1(x) \\ c_1(x) \\s_2(x) \\ c_2(x) \\ \vdots \end{pmatrix},\qquad \bfq(t)=\begin{pmatrix} \beta_1(x) \\ \alpha_1(x) \\ \beta_2(x) \\ \alpha_2(x) \\ \vdots \end{pmatrix},\qquad \bfp(t)=\begin{pmatrix} \eta_1(x) \\ \xi_1(x) \\ \eta_2(x) \\ \gamma_2(x) \\ \vdots \end{pmatrix}.
\end{equation}
Moreover, by setting $I_2$ the identity matrix of dimension 2, $J_2$ the skew-symmetric and orthogonal matrix defined in (\ref{J2dH}), and the infinite matrix
\begin{equation}\label{D}
D=\frac{2\pi}{b-a}\begin{pmatrix} ~1\\  & 2\\ &  & 3\\ &&&\ddots~ \end{pmatrix},
\end{equation}
the required partial derivatives of $u(x,t)$ and $v(x,t)$ can be easily computed as follows:
\begin{eqnarray}\nonumber
v_t(x,t) &=& \bfw(x)^\top \dot\bfp(t),\qquad v_x(x,t) = \left[(D\otimes J_2)\bfw(x)\right]^\top\bfp(t),\\
\label{der}
u_t(x,t)&=&\bfw(x)^\top \dot{\bfq}(t), \qquad u_x(x,t) = \left[(D\otimes J_2)\bfw(x)\right]^\top\bfq(t),\\
\nonumber
u_{xx}(x,t) &=& \left[(D\otimes J_2)^2\bfw(x)\right]^\top\bfq(t), \qquad
u_{xxx}(x,t) = \left[(D\otimes J_2)^3\bfw(x)\right]^\top\bfq(t),
\end{eqnarray}
due to the fact that
\begin{equation}\label{w1}
\bfw'(x) = (D\otimes J_2)\bfw(x).
\end{equation}
Consequently, by also considering that
\begin{equation}\label{J2pw}
J_2^\top = -J_2, \qquad J_2^2 = -I_2,  \qquad \int_a^b \bfw(x) \dd x = \bfzero , \qquad \int_a^b \bfw(x)\bfw(x)^\top \dd x = I,
\end{equation}
with $\bf 0$ the zero vector and $I$ the identity operator, the following result can be proved. 

\begin{theo}\label{fform} System (\ref{uvform}) can be cast in Hamiltonian form as
\begin{equation}\label{qpform}
\dot{\begin{pmatrix} \bfq\\ \bfp\end{pmatrix}} = \begin{pmatrix} (D\otimes J_2^\top)\bfp \\ (D\otimes J_2^\top) \left[
(D^2\otimes I_2)\bfq + \int_a^b \bfw(x)(\hat u_0+\bfw(x)^\top\bfq)^2\dd x\right] \end{pmatrix} 
\equiv \left(|J_2|\otimes D\otimes J_2^\top\right) \nabla H(\bfq,\bfp), 
\end{equation}
with Hamiltonian
\begin{equation}\label{Hfun1}
H(\bfq,\bfp) = \frac{1}2\left[ \bfp^\top\bfp + \bfq^\top(D^2\otimes I_2)\bfq +\frac{2}3\int_a^b (\hat u_0+\bfw(x)^\top\bfq)^3\dd x\right].
\end{equation}
This latter, in turn, is equivalent, up to a constant, to the functional (\ref{Hfun}), via the transformations (\ref{uv})--(\ref{der}).
\end{theo}
\proof The proof of (\ref{qpform}) follows by considering that, from (\ref{uv})--(\ref{J2pw}), one has, for the first equation,
\begin{eqnarray*}
\dot\bfq &=&\int_a^b \bfw(x) u_t\dd x = \int_a^b \bfw(x) v_x\dd x = \int_a^b \bfw(x) \left( (D\otimes J_2)\bfw(x)\right)^\top \bfp\,\dd x \\
&=& \underbrace{\int_a^b\bfw(x)\bfw(x)^\top \dd x}_{=\,I} \, (D\otimes J_2^\top)\bfp \,=\,(D\otimes J_2^\top)\bfp.
\end{eqnarray*}
Similarly, for second equation, by considering that, from (\ref{w1}) one has, by integrating by parts,
$$\int_a^b \bfw(x) (u^2)_x\dd x = \underbrace{\left[ \bfw(x) u^2\right]_{x=a}^{x=b}}_{=\,\bfzero} - \int_a^b (D\otimes J_2) \bfw(x) u^2\dd x =  (D\otimes J_2^\top)\int_a^b \bfw(x) u^2\dd x,$$
one obtains:
\begin{eqnarray*}
\dot\bfp &=&\int_a^b \bfw(x) v_t\dd x = \int_a^b \bfw(x) [ -u_{xxx} + (u^2)_x]\dd x = \underbrace{\int_a^b \bfw(x) \bfw(x)^\top \dd x}_{=\,I} (D^3\otimes J_2^\top)\bfq \\
&& +(D\otimes J_2^\top)\int_a^b \bfw(x) u^2\dd x = (D\otimes J_2^\top)\left[ (D^2\otimes I_2)\bfq + \int_a^b \bfw(x) (\hat u_0+\bfw(x)^\top \bfq)^2\dd x\right].
\end{eqnarray*}
The equivalence of (\ref{Hfun1}) with (\ref{Hfun}), up to a constant,  is explained below
\begin{eqnarray*}
\H[u,v] &=& \frac{1}2\int_a^b v^2 + u_x^2 +\frac{2}3u^3\dd x ~=~ \frac{1}2\int_a^b\Big\{ \left[ \hat v_0+\bfp^\top \bfw(x)\right] \left[ \hat v_0+\bfw(x)^\top \bfp\right] + \\[2mm]
&&\bfq^\top (D \otimes J_2)\bfw(x)\bfw(x)^\top (D \otimes J_2^\top)\bfq 
+ \frac{2}3(\hat u_0+\bfw(x)^\top\bfq)^3\Big\}\dd x\\
%\end{eqnarray*}\begin{eqnarray*}
&=& \frac{1}2\left[  \hat v_0^2\int_a^b\dd x\,+\,2\hat v_0\underbrace{\int_a^b\bfw(x)^\top\dd x}_{=\,\bfzero ^\top}\, \bfp\, +\, \bfp^\top \underbrace{\int_a^b\bfw(x)\bfw(x)^\top\dd x}_{=\,I}\, \bfp  \right.\\
&&+\left.\bfq^\top (D \otimes J_2)\underbrace{\int_a^b\bfw(x)\bfw(x)^\top\dd x}_{=\,I}\, (D \otimes J_2^\top)\bfq + \frac{2}3\int_a^b(\hat u_0+\bfw(x)^\top\bfq)^3\dd x\right]\\
&=& \frac{1}2\left[  (b-a)\hat v_0^2+\bfp^\top \bfp + 
\bfq^\top (D^2\otimes J_2 J_2^\top)\bfq + \frac{2}3\int_a^b(\hat u_0+\bfw(x)^\top\bfq)^3\dd x\right]\\
&=& \frac{1}2\left[  (b-a)\hat v_0^2+\bfp^\top \bfp + 
\bfq^\top (D \otimes I_2)\bfq + \frac{2}3\int_a^b(\hat u_0+\bfw(x)^\top\bfq)^3\dd x\right]\\
&=& \frac{b-a}2\hat v_0^2+H(\bfq,\bfp).\QED
\end{eqnarray*}

\medskip
Finally, by using similar arguments, the following result can be proved. 

\begin{theo}\label{Mth} The quadratic invariant (\ref{M}) is equivalent, up to a constant, to
\begin{equation}\label{Mqp} M(\bfq,\bfp) = \bfq^\top\bfp.\end{equation}
\end{theo}
\proof From (\ref{M}), (\ref{uv}), and (\ref{J2pw}), one has:
\begin{eqnarray*}
\M[u,v] &=& \int_a^b uv\dd x = \int_a^b (\hat u_0+\bfq^\top\bfw(x))(\hat v_0+\bfw(x)^\top \bfp)\dd x\\
&=&  \int_a^b \hat u_0\hat v_0\dd x + (\hat u_0\bfp+\hat v_0\bfq)^\top\underbrace{\int_a^b \bfw(x)\dd x}_{=\,\bfzero} \,+\,\bfq^\top \int_a^b \bfw(x)\bfw(x)^\top\dd x \, \bfp \\[2mm]
&=& (b-a)\hat u_0\hat v_0 +  \bfq^\top\bfp ~=~ (b-a)\hat u_0\hat v_0 + M(\bfq,\bfp).\,\QED
\end{eqnarray*}

\medskip
As is clear, in order to obtain a computational method, the infinite series in (\ref{uv}) have to be truncated at a convenient index $N$. In so doing, the infinite vectors and matrices in (\ref{wqp})-(\ref{D}) become of dimension $N$, i.e., respectively, 
\begin{equation}\label{wqpDN}
\bfw(x)=\begin{pmatrix} s_1(x) \\ c_1(x) \\ \vdots \\ s_N(x) \\ c_N(x) \end{pmatrix},\quad \bfq(t)=\begin{pmatrix} \beta_1(x) \\ \alpha_1(x) \\ \vdots\\ \beta_N(x) \\ \alpha_N(x) \end{pmatrix},\quad \bfp(t)=\begin{pmatrix} \eta_1(x) \\ \xi_1(x) \\ \vdots\\ \eta_N(x) \\ \gamma_N(x)\end{pmatrix}, \quad
D=\frac{2\pi}{b-a}\begin{pmatrix} 1\\  &\ddots \\ &  & N\end{pmatrix},
\end{equation}
Consequenly, (\ref{der}) continue formally to hold, even though now the truncated approximations to $u$ and $v$ do not satisfy the equations (\ref{uvform})  anymore. Nevertheless, in the spirit of Galerkin methods, by imposing the residual be orthogonal to the functional space spanned by the entries of (the truncated version of) $\bfw(x)$, the results of Theorems~\ref{fform} and \ref{Mth} continue formally to hold, with the only difference that now the truncated versions of $H(\bfq,\bfp)$ and $M(\bfq,\bfp)$ do not coincide, up to a constant, with the functionals (\ref{Hfun}) and (\ref{M}), respectively. Nevertheless, it is known that, upon regularity assumptions on  $u$ and $v$, the truncated series
$(\hat u_0+\bfw(x)^\top\bfq)$ and $(\hat v_0+\bfw(x)^\top\bfp)$ converge more than exponentially to them, as well as the truncated version of $H(\bfq,\bfp)$ and $M(\bfq,\bfp)$ to the corresponding functionals, as $N\rightarrow\infty$. This phenomenon is usually referred to as to {\em spectral accuracy} (see, e.g., \cite{Tref}).

For completeness, we mention that, in order to obtain a fully semi-discrete problem, the integrals appearing in (\ref{qpform}) and (\ref{Hfun1}) need to be evaluated. In the present case, since the Hamiltonian is a polynomial of degree 3, this can be done exactly (see, e.g., \cite[Theorem~7]{BFCI2015}) by using a composite trapezoidal rule based at the evenly spaced points:
\begin{equation}\label{xim}
x_i = a+i\frac{b-a}m, \qquad i=0,1,\dots,m,
\end{equation}
with $m=2N+1$ (for (\ref{qpform})) and $m=3N+1$ (for (\ref{Hfun1})), respectively.

\section{Hamiltonian Boundary Value Methods}\label{HBVM}
In this section, we recall the main facts about Hamiltonian Boundary Value Methods (HBVMs), which constitute a class of energy-conserving Runge-Kutta methods for Hamiltonian problems. Moreover, we study their efficient implementation for solving problem (\ref{qpform})-(\ref{wqpDN}). HBVMs have been investigated in a series of papers \cite{BIT2010,BIT2012,BIT2015} (see also the monograph \cite{LIMbook2016} and the recent review paper \cite{BI2018}) for solving Hamiltonian problems, and have been developed in a series of directions (see, e.g., \cite{BCMR2012,BGIW2018,BI2012,BIT2012-1}), including Hamiltonian boundary value problems \cite{ABI2015}. More recently, they have been successfully used to solve Hamiltonian PDEs \cite{BBFCI2018,BFCI2015,BGS2018,LIMbook2016,BI2018,BZL2018}, and this paper belongs to this last field of investigation.

In more detail, for all $k\ge s$, the HBVM$(k,s)$ method is the $k$-stage Runge-Kutta method with Butcher tableau given by
\begin{equation}\label{butab}
\begin{array}{c|c}
\bfc & \I_s\P_s^\top\Omega \\ \hline \\[-3mm] &\bfb^\top
\end{array}\,,
\end{equation}
where, by setting $\{P_j\}$ the Legendre polynomials shifted and scaled in order to be orthonormal on the interval $[0,1]$,
\begin{eqnarray}\nonumber
\bfc &=& (c_1,\dots,c_k)^\top, \quad \bfb~=~(b_1,\dots,b_k)^\top, \qquad \Omega = \begin{pmatrix} b_1\\ &\ddots\\ &&b_k\end{pmatrix},\\
\label{butab1}\\ \nonumber
\I_s &=& \begin{pmatrix} \int_0^{c_1} P_0(x)\dd x& \dots &\int_0^{c_1} P_{s-1}(x)\dd x\\
\vdots & &\vdots\\
\int_0^{c_k} P_0(x)\dd x& \dots &\int_0^{c_k} P_{s-1}(x)\dd x\end{pmatrix}, \quad
\P_s ~=~ \begin{pmatrix} P_0(c_1) &\dots &P_{s-1}(c_1)\\
\vdots & &\vdots\\ P_0(c_k) &\dots &P_{s-1}(c_k)\end{pmatrix},
\end{eqnarray}
with $(c_i,b_i)$ the Legendre abscissae and weights of the Gauss interpolatory quadrature formula of order $2k$ (i.e., $P_k(c_i)=0$, $i=1,\dots,k$).

When applied for solving the ODE-IVPs
\begin{equation}\label{ivp}
\dot y = f(y), \qquad y(0)=y_0\in\RR^m,
\end{equation}
with a stepsize $h$, the Runge-Kutta method (\ref{butab})-(\ref{butab1}) implicitly defines a polynomial approximation $\sigma\in\Pi_s$ such that 
\begin{equation}\label{lim}
\sigma(0)=y_0, \qquad \sigma(h)=:y_1\approx y(h), 
\end{equation}
providing, in case of Hamiltonian problems, relevant conservation properties, as is specified by the following theorem  \cite{LIMbook2016,BIT2012,BIT2015,BI2018}.

\begin{theo}\label{HBVMth} For all $k\ge s$, a HBVM$(k,s)$ method used with stepsize $h$:
\begin{itemize}
\item is symmetric and  \,$y_1-y(h)=O(h^{2s+1})$;
\item when $k=s$, it coincides with the symplectic $s$-stage Gauss method.
\end{itemize}
Moreover, when solving an Hamiltonian system,  i.e., in (\ref{ivp}) $f(y)=J\nabla H(y)$ with $J^\top=-J$:
\begin{itemize}
\item it is energy-conserving when the Hamiltonian $H$ is a polynomial and \,$\deg H \le 2k/s$\,;
\item conversely, one has \,$H(y_1)-H(y_0)=O(h^{2k+1})$.
\end{itemize}
\end{theo}
\begin{rem} From the last two points in Theorem~\ref{HBVMth}, one has that, by 
choosing $k$ large enough, either an exact energy-conservation can be gained,
in the polynomial case, or a {\em ``practical''} energy-conservation can be obtained in the general case.
In fact, in the latter case, it is enough that the energy error falls within the round-off
error level.
\end{rem}

As a consequence, one has the following conservation result.

\begin{cor}\label{HBVMGBE} 
For all $k\ge \frac{3}2s$, the HBVM$(k,s)$ method is energy-conserving and of order $2s$, when used for solving the Hamiltonian problem (\ref{qpform})-(\ref{Hfun1}).\end{cor}

\proof In fact, in such a case, $\deg H = 3\le 2k/s$, for all $k\ge \frac{3}2s$.\,\QED\medskip

We now sketch the efficient implementation of HBVMs, in view of their application for solving the semi-discrete problem (\ref{qpform})-(\ref{wqpDN}). 
To begin with, one of the main features of a HBVM$(k,s)$ method is that the discrete problem generated by the application of the method has (block) dimension $s$, independently of $k$. This feature, in turn, allows the use of possibly much larger values of $k$, w.r.t. $s$. The key point for this \cite{BIT2011} stems from the fact that the polynomial approximation (\ref{lim}) has degree $s$, and the discrete problem can be cast in terms of its unknown coefficients. In more details, the equation for the $k$ stages $Y_1,\dots,Y_k$ of the Runge-Kutta method (\ref{butab}) solving (\ref{ivp}) can be written as
\begin{equation}\label{stages}
Y = \bfe\otimes y_0+h\I_s\P_s^\top\Omega\otimes I_m f(Y),
\end{equation}
having set $$Y = \begin{pmatrix} Y_1\\ \vdots \\ Y_k\end{pmatrix}, \qquad f(Y) = \begin{pmatrix} f(Y_1)\\ \vdots \\ f(Y_k)\end{pmatrix}, \qquad \bfe = \begin{pmatrix} 1\\ \vdots\\ 1\end{pmatrix}\in\RR^k,$$
with the new approximation given by
\begin{equation}\label{y1}
y_1 = y_0 + h\sum_{i=1}^k b_i f(Y_i).
\end{equation}
By defining the block vector of dimension $s$:
\begin{equation}\label{gamma}
\bfgamma = \begin{pmatrix} \gamma_0\\ \vdots \\ \gamma_{s-1}\end{pmatrix} := \P_s^\top\Omega\otimes I_m f(Y),
\end{equation}
one then obtains that (\ref{stages}) can be re-written as:
\begin{equation}\label{Y}
Y = \bfe\otimes y_0+h\I_s\otimes I_m \bfgamma,
\end{equation}
which, substituted at the right-hand side in (\ref{gamma}), provides us with the following equation,
\begin{equation}\label{gamma1}
\bfgamma = \P_s^\top\Omega\otimes I_m f\left(\bfe\otimes y_0+h\I_s\otimes I_m \bfgamma\right),
\end{equation}
having block dimension $s$. It can be easily seen that:
\begin{itemize}
\item the polynomial $\sigma$ mentioned in (\ref{lim}) is given by
\begin{equation}\label{sigma}
\sigma(ch) = y_0 + h\sum_{j=0}^{s-1}\int_0^c P_j(x)\dd x \gamma_j, \qquad c\in[0,1];
\end{equation}
\item the stages of the $k$-stage HBVM$(k,s)$ method are defined by
$$Y_i = \sigma(c_ih), \qquad i=1,\dots,k,$$
\item the new approximation (\ref{y1}) is given by, setting $c=1$ in (\ref{sigma}):
\begin{equation}\label{y1new}
y_1 = y_0 + h\gamma_0.
\end{equation}
\end{itemize}
Consequently, in order to implement the step of an HBVM$(k,s)$ method, one needs to solve the discrete problem (\ref{gamma1}), i.e., the equation
\begin{equation}\label{F}
F(\bfgamma) := \bfgamma - \P_s^\top\Omega\otimes I_m f\left(\bfe\otimes y_0+h\I_s\otimes I_m \bfgamma\right) = \bfzero.
\end{equation}
This equation, which could in principle be solved by using a straightforward fixed-point iteration, 
$$\bfgamma^{\ell+1} = \P_s^\top\Omega\otimes I_m f\left(\bfe\otimes y_0+h\I_s\otimes I_m \bfgamma^\ell\right), \qquad \ell=0,1,\dots,$$
actually requires, in the case of problem (\ref{qpform})-(\ref{wqpDN}), the use of a Newton-type iteration, in order to avoid the use of very small stepsizes (indeed, of the order of
$\|D\|^{-3}\propto N^{-3}$, which becomes very small when large values of $N$ are considered). For this purpose, let us consider the simplified Newton iteration
for solving (\ref{F}) which, by considering that (see, e.g., \cite{BIT2011,LIMbook2016})
\begin{equation}\label{Xs}
\P_s^\top\Omega\I_s = X_s := \begin{pmatrix} \xi_0 & -\xi_1\\ \xi_1 & 0 & \ddots \\ &\ddots &\ddots &-\xi_{s-1}\\ & & \xi_{s-1} & 0\end{pmatrix}, \qquad \xi_i = \left(2\sqrt{|4i^2-1|}\right)^{-1}, \quad i=0,\dots,s-1,
\end{equation}
formally reads
\begin{equation}\label{sNit}
\bfgamma^{\ell+1} = \bfgamma^\ell - \left[ I-hX_s\otimes f'(y_0)\right]^{-1} F(\bfgamma^\ell), \qquad \ell=0,1,2,\dots,
\end{equation}
with $f'$ the Jacobian of the function $f$ in (\ref{ivp}). This iteration, in turn, requires the factorization of the matrix
\begin{equation}\label{mat}
\left[ I-hX_s\otimes f'(y_0)\right]
\end{equation}
having dimension $s$ times larger than that of $f'(y_0)$. It can be proved that this iteration can be conveniently replaced by a corresponding {\em blended iteration} \cite{Br2000,BrMa2002,BrMa2009,BIT2011,LIMbook2016}
which, having set\,\footnote{As is usual, $\sigma(X_s)$ denotes the spectrum of matrix $X_s$.}
\begin{equation}\label{Sigma}
\Sigma = I_m - h\rho_s f'(y_0), \qquad\mbox{with}\qquad \rho_s =\min_{\lambda\in\sigma(X_s)} |\lambda|,
\end{equation}
reads:
\begin{eqnarray*}
&& \bfeta^\ell = -F(\bfgamma^\ell), \qquad \bfeta_1^\ell = (\rho_s X_s^{-1}\otimes I_m)\bfeta^\ell,\\[2mm]
&&\bfgamma^{\ell+1} = \bfgamma^\ell - I_s\otimes \Sigma^{-1}\left[ \bfeta_1^\ell + I_s\otimes \Sigma^{-1}\left( \bfeta^\ell-\bfeta_1^\ell\right)\right], \qquad \ell=0,1,\dots.
\end{eqnarray*}
Consequently, only the factorization of matrix $\Sigma$ in (\ref{Sigma}) is needed, {\em independently of $s$}. This, in turn, allows the use of relatively large values of $s$.\footnote{This feature will be very important, as we are going to see in the sequel.}
Moreover, in the case of the problem (\ref{qpform})-(\ref{wqpDN}), one has the further simplification that the Jacobian of the right-hand side can be conveniently approximated by the linear part alone, i.e.,
$$%\begin{equation}\label{apprJ}
\begin{pmatrix} & (D\otimes J_2^\top)\\ (D^3\otimes J_2^\top) \end{pmatrix},
$$%\end{equation}
so that matrix $\Sigma$ in (\ref{Sigma}) becomes, by setting hereafter $I\in\RR^{N\times N}$ the identity matrix and $h$ the used time step,
\begin{equation}\label{Sigma1}
\Sigma = \begin{pmatrix} I\otimes I_2& -\tau (D\otimes J_2^\top)\\ -\tau (D^3\otimes J_2^\top) &I\otimes I_2 \end{pmatrix},  \qquad\mbox{where}\qquad \tau = \rho_s h.
\end{equation}
Consequently, its inverse has to be computed only once for all.

It must be stressed that the particular structure of the Jacobian matrix and, therefore, of matrix $\Sigma$, is of paramount importance to derive an efficient numerical method, based on the use of HBVMs.
Clearly, this structure is strictly related to the equation at hand, and this fact has been investigated for many Hamiltonian PDEs \cite{BFCI2015,BBFCI2018,BZL2018}. In particular, for the ``good'' Boussinesq equation, 
the following result holds true.

\begin{theo}\label{Sinv}
With reference to matrix $\Sigma$ defined in (\ref{Sigma1}), one has:
\begin{equation}
\Sigma^{-1} =  \mat D_1\otimes I_2& \tau (D_1D\otimes J_2^\top)\\ \tau (D_1D^3\otimes J_2^\top) &D_1\otimes I_2\endmat, \qquad\mbox{with}\qquad
D_1 = (I+\tau^2D^4)^{-1}.   \label{Sinv1}
\end{equation}
\end{theo}
\proof
In fact, one has: 
$$\mat I\otimes I_2 \\ &I\otimes J_2^\top\endmat \Sigma \mat I\otimes I_2 \\ &I\otimes J_2\endmat = \mat I &-\tau D\\ \tau D^3 &I\endmat \otimes I_2$$
and
$$\mat I &-\tau D\\ \tau D^3 &I\endmat^{-1} = (I_2\otimes D_1) \mat I & \tau D\\ -\tau D^3 & I\endmat,$$
with $D_1$  defined as in (\ref{Sinv1}). Consequently,
\begin{eqnarray*}
\Sigma^{-1} &=&  \mat I\otimes I_2 \\ &I\otimes J_2\endmat 
(I_2\otimes D_1\otimes I_2) \left[\mat I & \tau D\\ -\tau D^3 & I\endmat\otimes I_2\right] \mat I\otimes I_2 \\ &I\otimes J_2^\top\endmat\\[2mm]
&=&\left(I_2\otimes D_1\otimes I_2\right) \mat I\otimes I_2& \tau (D\otimes J_2^\top)\\ \tau (D^3\otimes J_2^\top) &I\otimes I_2\endmat
=\mat D_1\otimes I_2& \tau (D_1D\otimes J_2^\top)\\ \tau (D_1D^3\otimes J_2^\top) &D_1\otimes I_2\endmat.
\,\QED
\end{eqnarray*}
\medskip

As a consequence of the previous theorem, one has that matrix $\Sigma^{-1}$ in (\ref{Sinv1}) can be computed with a cost which is linear in the dimension of the problem, since it has blocks with a diagonal structure.
Moreover, it can be conveniently stored by using 3 vectors of dimension $N$ (containing the diagonal entries of $D_1$, $\tau D_1D$, and $\tau D_1D^3$).

\subsection{Spectral HBVMs}\label{shbvms}
The previous blended implementation of HBVMs is particularly interesting, since it allows the use of relatively large values of $s$. This, in turn, allows to use HBVMs as spectral methods in time \cite{BMR2018,BIMR2018}, so that one obtains, for the used finite precision arithmetic, the maximum possible accuracy compatible with the considered timestep. We here sketch the use of HBVMs as spectral methods (which we shall refer to as {\em spectral HBVMs or, in short, SHBVMs}): further details can be found in the previous references \cite{BMR2018,BIMR2018}. 

To begin with, let us consider the problem (\ref{ivp}) which, by considering the time interval $[0,h]$ and expanding the right-hand side along the Legendre basis, can be rewritten as
$$\dot y(ch) = \sum_{j\ge0} P_j(c) \gamma_j(y), \qquad c\in[0,1], \qquad \gamma_j(y)=\int_0^1P_j(c)f(y(ch))\dd c, \qquad j=0,1,\dots.$$
Then, integrating term by term and imposing that $y(0)=y_0$, one obtains that the solution is formally given by:
$$y(ch) = y_0 + h\sum_{j\ge0} \int_0^cP_j(x)\dd x\, \gamma_j(y), \qquad c\in[0,1].$$
Furthermore, for a suitably regular function $f$ one has that 
\begin{equation}\label{vaa0}
\gamma_j(y)\rightarrow 0, \qquad\mbox{as}\qquad j\rightarrow\infty.
\end{equation}
Consequently, $y(ch)$ can be approximated within machine accuracy by the polynomial (\ref{sigma}), provided that the quadrature $(c_i,b_i)$ is accurate enough (i.e., $k$ is large enough), and $s$ is the first index such that
\begin{equation}\label{tol}
\|\gamma_{s-1}\|\le tol ~\cdot\max_{j=0,\dots,s-1} \|\gamma_j\|,
\end{equation}
with $tol\sim u$, being $u$ the machine epsilon of the considered finite precision arithmetic.

In the practice, however, the numerical evaluation of the coefficients $\gamma_j\approx \gamma_j(y)$, makes them ``stagnate'' (in norm) around a small value, rather than tending to 0, according to (\ref{vaa0}): in such a case, the tolerance $tol$ in (\ref{tol}) is more conveniently chosen in order to avoid using the coefficients with a stagnating norm, since this means that they are not reliably computed. This criterion will always be used in the sequel, for implementing SHBVMs.

\section{Numerical tests}\label{NT}

In this section, we report a few numerical tests concerning the numerical solution of problem (\ref{qpform})-(\ref{wqpDN}) with initial conditions given by (see (\ref{uv0})):
\begin{equation}\label{q0p0}
\bfq(0) = \int_a^b \bfw(x)u_0(x)\dd x, \qquad \bfp(0) = \int_a^b \bfw(x)v_0(x)\dd x.
\end{equation}
In particular, we compare the following methods:
\begin{itemize}
\item the symplectic $s$-stage Gauss methods (which we shall denote Gauss $s$), having order $2s$, for  $s=1,2$. Such methods are expected to conserve the quadratic invariant (\ref{Mqp}) but only approximately the Hamiltonian (\ref{Hfun1});
\item the energy-conserving HBVM(2,1) and HBVM(3,2) methods, having respectively order $2$ and $4$. Such method conserve the Hamiltonian (\ref{Hfun1}) but only approximately the quadratic invariant (\ref{Mqp});
\item the spectral HBVM (SHBVM) method, using a value of $s$ and $k=\lceil1.5s\rceil$ large enough so that the maximum possible accuracy is gained. As sketched in Section~\ref{shbvms}, the value of $s$ is obtained by appropriately choosing the tolerance $tol$ in (\ref{tol}). Such methods are expected to conserve both the Hamiltonian (\ref{Hfun1}) and the momentum (\ref{Mqp}), as well as to provide a solution error within the round-off error level.
\end{itemize}
For each considered problem, we compare the above methods in terms of:
\begin{itemize}
\item the maximum solution error $e_u$;
\item the maximum Hamiltonian error $e_H$;
\item the maximum momentum error $e_M$;
\item the execution time (in sec);
\item moreover, when appropriate, we also estimate the numerical rate of convergence.
\end{itemize}
All numerical tests have been performed on a 3.1GHz quad-core Intel i7 computer with 16GB of memory, running Matlab 2017b. Moreover, the same Matlab code implements all the above methods, so that the comparisons are fair. In all cases, the blended iteration previously described has been used.

\subsection*{Solitary wave}
\begin{figure}[t]
\centerline{\includegraphics[height=9cm,width=12cm]{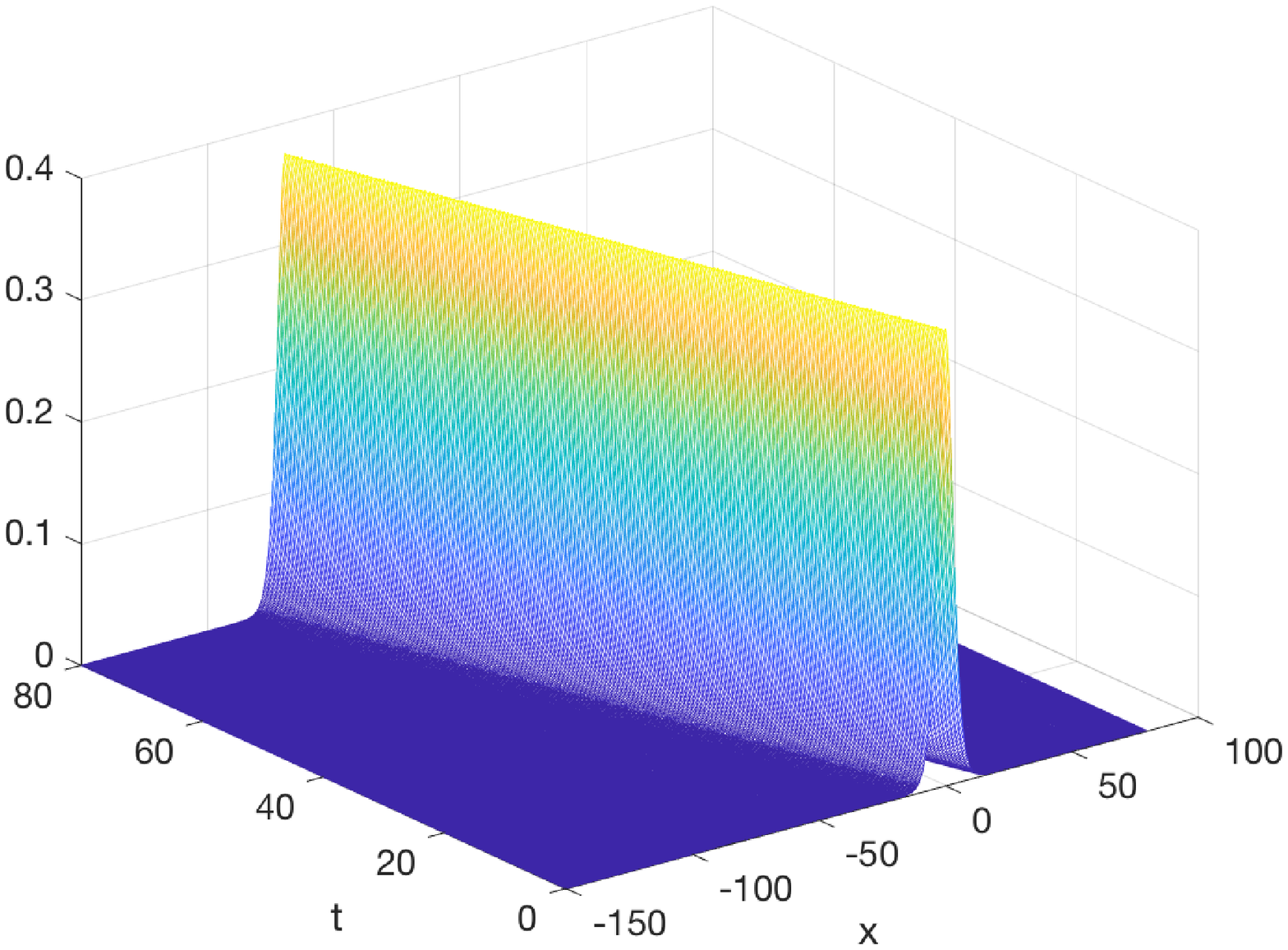}}
\caption{Plot of $\frac{1}2-u(x,t)$ for problem (\ref{uvform}) and (\ref{swsol})-(\ref{swuv0}).}
\label{swfig}
\end{figure}

Let us at first consider the  solitary wave solution \cite{Chen2} of (\ref{uvform}) given, by taking into account (\ref{u2w}), by 
\begin{equation}\label{swsol}
u(x,t) = \frac{1}2-A\cdot \sech^2\left( \sqrt{\frac{A}6}(x+ct-\xi_0)\right), \qquad  c=\pm\sqrt{1-\frac{2}3A}. %A = \frac{3}2P^2, \qquad c=\pm\sqrt{1-P^2}.
\end{equation}
Consequently, the initial conditions (\ref{uv0}) at $t=0$ are given by:
\begin{equation}\label{swuv0}
u_0(x) = \frac{1}2-A\cdot \sech^2\left( \sqrt{\frac{A}6}(x-\xi_0)\right), \qquad v_0(x) = c\left(u_0(x)-\frac{1}2\right).
\end{equation}
We consider the values $\xi_0=0$, $A=3/8$, and the positive value of $c$. We integrate in time until $T=80$, so that if we consider the space interval $[-120,80]$ both $u$ and $v$ can be assumed to be approximately periodic.\footnote{Actually, ``exactly'' periodic, when using the double precision IEEE.} The expansions (\ref{uv}) have been truncated at $N=300$, providing spectral accuracy in space. As matter of fact, the spatial semi-discretization error, measured on the initial conditions (see (\ref{q0p0})), which is defined as
\begin{equation}\label{e0}
e_0 := \max\left\{ \| u_0(x)-\hat u_0 -\bfw(x)^\top\bfq(0)\|,\, \| v_0(x)-\hat v_0-\bfw(x)^\top\bfp(0)\| \right\},
\end{equation}
is $5.06\times 10^{-14}$. The solution (\ref{swsol}) of the problem is depicted in Figure~\ref{swfig}, whereas in Table~\ref{swtab} we list the obtained numerical results, as explained above, by using a timestep $h = 80/n$. For the SHBVM method, we used a tolerance $tol\sim 10^{-11}$ in (\ref{tol}), providing $s=10$ (and, therefore, $k=15$). From the results reported in Table~\ref{swtab}, one infers that the latter method (SHBVM) is the most effective one among those considered, able to numerically conserve all the invariants, while providing a negligible solution error, with a very small execution time.

\begin{table}[p]
\caption{Numerical results with  timestep $h=80/n$ for the solitary wave problem (\ref{uvform}) and (\ref{swsol})-(\ref{swuv0}).}
\label{swtab}

\smallskip
\centerline{\begin{tabular}{|r|r|c|c|c|c|c|c|}
\hline
\multicolumn{8}{|c|}{Gauss 1}\\
\hline
$n$ & time & $e_u$ & rate & $e_H$ & rate & $e_M$ & ---\\ 
\hline
8000 & 34.3 & 3.87e-06 &  --- & 2.80e-12&  ---  & 2.78e-15 &   \\ 
9600 & 40.8 & 2.69e-06 &  2.0 & 1.37e-12&  3.9  & 2.55e-15 &    \\ 
11200 & 47.8 & 1.97e-06 &  2.0 & 7.62e-13&  3.8  & 2.78e-15 &   \\ 
12800 & 54.4 & 1.51e-06 &  2.0 & 4.62e-13&  3.8  & 2.22e-15 &    \\ 
14400 & 61.6 & 1.19e-06 &  2.0 & 3.04e-13&  3.6  & 3.77e-15 &   \\ 
16000 & 76.6 & 9.67e-07 &  2.0 & 2.10e-13&  3.5  & 4.22e-15 &   \\ 
\hline
\hline
\multicolumn{8}{|c|}{Gauss 2}\\
\hline
$n$ & time & $e_u$ & rate & $e_H$ & rate & $e_M$ & ---\\ 
\hline
1600 & 19.0 & 1.01e-09 & --- & 3.73e-14&  ---  & 1.67e-15 &    \\ 
2400 & 28.3 & 1.99e-10 &  4.0 & 3.55e-14&  **  & 1.67e-15 &    \\ 
3200 & 33.3 & 6.30e-11 &  4.0 & 3.73e-14& **  & 1.78e-15 &   \\ 
4000 & 41.8 & 2.58e-11 &  4.0 & 4.44e-14& **  & 1.89e-15 &   \\ 
\hline
\hline
\multicolumn{8}{|c|}{HBVM(2,1)}\\
\hline
$n$ & time & $e_u$ & rate & $e_H$ & --- & $e_M$ & rate\\ 
\hline
8000 & 54.3 & 3.97e-06 &  --- & 1.24e-14&     & 2.97e-12 &  --- \\ 
9600 & 66.0 & 2.76e-06 &  2.0 & 1.42e-14&    & 1.45e-12 &  3.9 \\ 
11200 & 76.8 & 2.03e-06 &  2.0 & 1.24e-14&     & 8.01e-13 &  3.9 \\ 
12800 & 87.7 & 1.55e-06 &  2.0 & 1.24e-14&     & 4.89e-13 &  3.7 \\ 
14400 & 98.6 & 1.23e-06 &  2.0 & 1.24e-14&     & 3.17e-13 &  3.7 \\ 
16000 & 109.1 & 9.93e-07 &  2.0 & 1.42e-14&    & 2.23e-13 &  3.4 \\ 
\hline
\multicolumn{8}{|c|}{HBVM(3,2)}\\
\hline
$n$ & time & $e_u$ & rate & $e_H$ & --- & $e_M$ & rate\\ 
\hline
1600 & 24.0 & 9.96e-10 &  --- & 1.07e-14&     & 3.97e-14 &  --- \\ 
2400 & 32.1 & 1.97e-10 &  4.0 & 1.07e-14&     & 3.96e-14 &  ** \\ 
3200 & 37.8 & 6.22e-11 &  4.0 & 1.42e-14&    & 4.25e-14 & ** \\ 
4000 & 47.3 & 2.55e-11 &  4.0 & 1.24e-14&    & 4.02e-14 &  ** \\ 
\hline
\hline
\multicolumn{8}{|c|}{SHBVM $(k=15,s=10)$}\\
\hline
$n$ & time & $e_u$ & --- & $e_H$ & --- & $e_M$ & ---\\ 
\hline
80  & 9.1   & 4.70e-14 & & 8.88e-15 & & 3.29e-14 & \\
\hline
\end{tabular}}
\end{table}

\subsection*{Spread of two solitary waves}
\begin{figure}[t]
\centerline{\includegraphics[height=9cm,width=12cm]{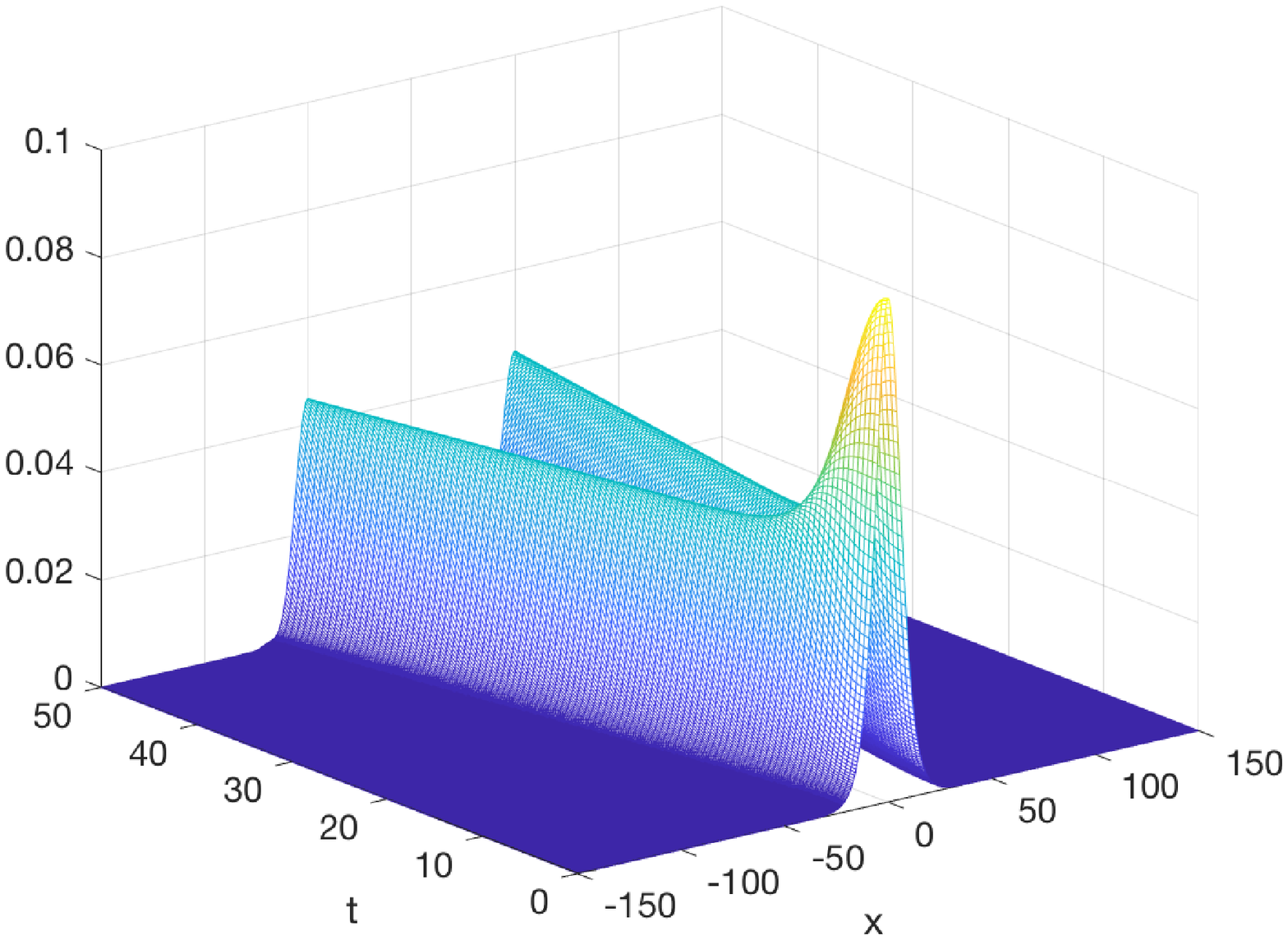}}
\caption{Plot of $\frac{1}2-u(x,t)$ for problem (\ref{uvform})-(\ref{dwuv0}).}
\label{dwfig}
\end{figure}

In general, the superposition of solitary waves as (\ref{swsol}) is no more a solution of (\ref{uvform}). Nevertheless, it provides an approximate solution configuration for that equation. As an example, the following initial conditions:
\begin{equation}\label{dwuv0}
u_0(x) = \frac{1}2-A\cdot \sech^2\left( \sqrt{\frac{A}6}x\right), \qquad v_0(x) \equiv 0,
\end{equation}
provide a single wave that, after a transient phase, approximately generates two solitary waves moving in opposite directions. We choose the parameters
$A = 3/32$, the space interval $[-150,150]$, and integrate until $T=50$. The expansions (\ref{uv}) have been truncated at $N=300$, providing a specrtal accuracy in space, with a spatial semi-discretization error (\ref{e0}) of $5.00\times 10^{-14}$. The corresponding solution is depicted in Figure~\ref{dwfig}. In Table~\ref{dwtab} we list the obtained numerical results, when using a timestep $h=50/n$. For the SHBVM method, we used a tolerance $tol\sim 10^{-10}$ in (\ref{tol}), again providing $s=10$ (and $k=15$).\footnote{The reference solution has been computed by using the SHBVM on a doubled time mesh.}  As in the previous example, this latter method turns out to be the most effective one, among those considered here, able to numerically conserve all the invariants and providing a negligible solution error, with a very small execution time. 

\begin{table}[p]
\caption{Numerical results with timestep $h=50/n$ for the spread of solitary waves problem (\ref{uvform})-(\ref{dwuv0}).}
\label{dwtab}

\smallskip
\centerline{\begin{tabular}{|r|r|c|c|c|c|c|c|}
\hline
\multicolumn{8}{|c|}{Gauss 1}\\
\hline
$n$ & time & $e_u$ & rate & $e_H$ & rate & $e_M$ & ---\\ 
\hline
5000 & 21.3 & 1.31e-07 &  --- & 1.66e-10&  0.0  & 4.23e-17 &    \\ 
6000 & 25.6 & 9.13e-08 &  2.0 & 1.15e-10&  2.0  & 4.08e-17 &    \\ 
7000 & 29.9 & 6.71e-08 &  2.0 & 8.48e-11&  2.0  & 3.96e-17 &    \\ 
8000 & 34.2 & 5.13e-08 &  2.0 & 6.49e-11&  2.0  & 4.04e-17 &   \\ 
9000 & 38.7 & 4.06e-08 &  2.0 & 5.13e-11&  2.0  & 3.92e-17 &    \\ 
10000 & 42.7 & 3.29e-08 &  2.0 & 4.16e-11&  2.0  & 3.94e-17 &   \\ 
\hline
\hline
\multicolumn{8}{|c|}{Gauss 2}\\
\hline
$n$ & time & $e_u$ & rate & $e_H$ & rate & $e_M$ & ---\\ 
\hline
1000 & 10.3 & 1.81e-11 &  --- & 7.11e-14&  ---  & 3.97e-17 &    \\ 
1500 & 15.4 & 3.57e-12 &  4.0 & 6.04e-14& **  & 4.00e-17 &   \\ 
2000 & 17.7 & 1.14e-12 &  4.0 & 6.57e-14& **  & 4.02e-17 &   \\ 
2500 & 22.0 & 4.80e-13 &  3.9 & 6.93e-14& **  & 4.10e-17 &  \\ 
\hline
\hline
\multicolumn{8}{|c|}{HBVM(2,1)}\\
\hline
5000 & 34.0 & 1.31e-07 &  --- & 3.20e-14&     & 3.88e-17 &  --- \\ 
6000 & 40.9 & 9.09e-08 &  2.0 & 3.73e-14&    & 4.05e-17 & ** \\ 
7000 & 47.7 & 6.68e-08 &  2.0 & 4.62e-14&    & 3.98e-17 &  ** \\ 
8000 & 54.3 & 5.11e-08 &  2.0 & 3.38e-14&     & 4.13e-17 & ** \\ 
9000 & 61.3 & 4.04e-08 &  2.0 & 3.91e-14&    & 4.10e-17 &  ** \\ 
10000 & 68.2 & 3.27e-08 &  2.0 & 3.55e-14&     & 4.00e-17 &  ** \\ 
\hline
$n$ & time & $e_u$ & rate & $e_H$ & --- & $e_M$ & rate\\ 
\hline
\hline
\multicolumn{8}{|c|}{HBVM(3,2)}\\
\hline
$n$ & time & $e_u$ & rate & $e_H$ & --- & $e_M$ & rate\\ 
\hline
1000 & 11.1 & 1.78e-11 &  --- & 3.38e-14&     & 3.97e-17 &  --- \\ 
1500 & 16.6 & 3.53e-12 &  4.0 & 3.91e-14&    & 4.11e-17 & ** \\ 
2000 & 19.1 & 1.13e-12 &  4.0 & 3.20e-14&     & 3.89e-17 &  ** \\ 
2500 & 23.9 & 4.80e-13 &  3.8 & 3.55e-14&    & 3.85e-17 &  ** \\ 
\hline
\hline
\multicolumn{8}{|c|}{SHBVM $(k=15,s=10)$}\\
\hline
$n$ & time & $e_u$ & --- & $e_H$ & --- & $e_M$ & ---\\ 
\hline
50   & 3.5   & 5.58e-14 & &1.60e-14  & & 4.07e-17 & \\
\hline
\end{tabular}}
\end{table}

\subsection*{Collision of two solitary waves}
\begin{figure}[t]
\centerline{\includegraphics[height=9cm,width=12cm]{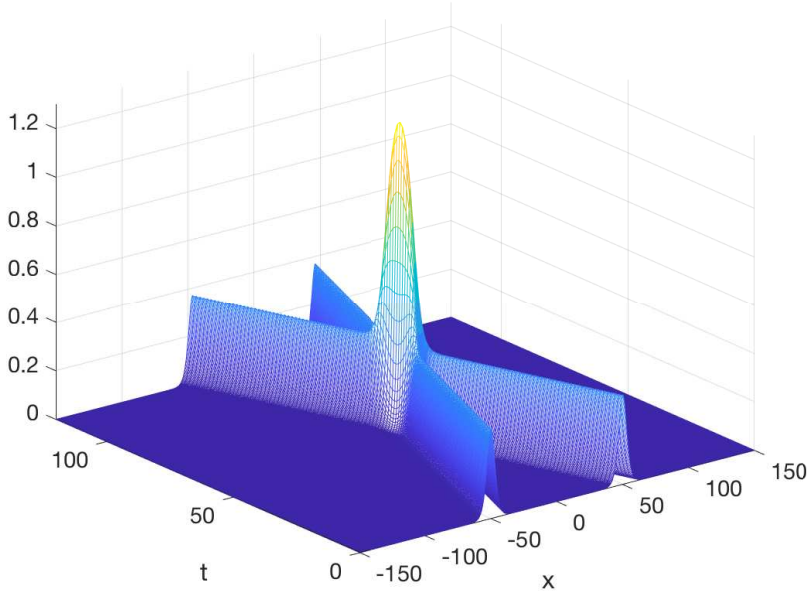}}
\caption{Plot of $\frac{1}2-u(x,t)$ for problem (\ref{uvform})-(\ref{cwuv0}).}
\label{cwfig}
\end{figure}

The last test problem we consider is provided by the following initial conditions,
\begin{eqnarray}\nonumber
u_0(x) &=& \frac{1}2-A\sech^2\left( \sqrt{\frac{A}6}(x-\xi_2)\right) -A\sech^2\left( \sqrt{\frac{A}6}(x-\xi_1)\right), \\ \label{cwuv0}\\ \nonumber
v_0(x) &=& c\left[A\sech^2\left( \sqrt{\frac{A}6}(x-\xi_2)\right) -A\sech^2\left( \sqrt{\frac{A}6}(x-\xi_1)\right)\right],
\end{eqnarray}
which, when choosing the parameters $A=0.369$, $c=\sqrt{1-\frac{2}3A}$, $\xi_2=-\xi_1=50$, provide two waves, which collide at about $t\approx60$.
We choose the space interval as $[-150,150]$ and integrate until $T=120$. The corresponding solution is depicted in Figure~\ref{cwfig}. The expansions (\ref{uv}) have been truncated at $N=300$, providing spectral accuracy in space, with a spatial semi-discretization error (\ref{e0}) of $5.01\times 10^{-14}$.  In Table~\ref{cwtab} we list the obtained numerical results, when using a timestep $h=120/n$. For the SHBVM method, we used a tolerance $tol\sim 10^{-11}$ in (\ref{tol}), providing $s=12$ (and, then, $k=18$).\footnote{The reference solution has been computed by using the SHBVM on a doubled time mesh.} As in the previous cases, this latter method turns out to be the most effective one, conserving all the invariant and with a negligible solution error, and a small execution time. 

\begin{figure}[t]
\centerline{\includegraphics[height=9cm,width=12cm]{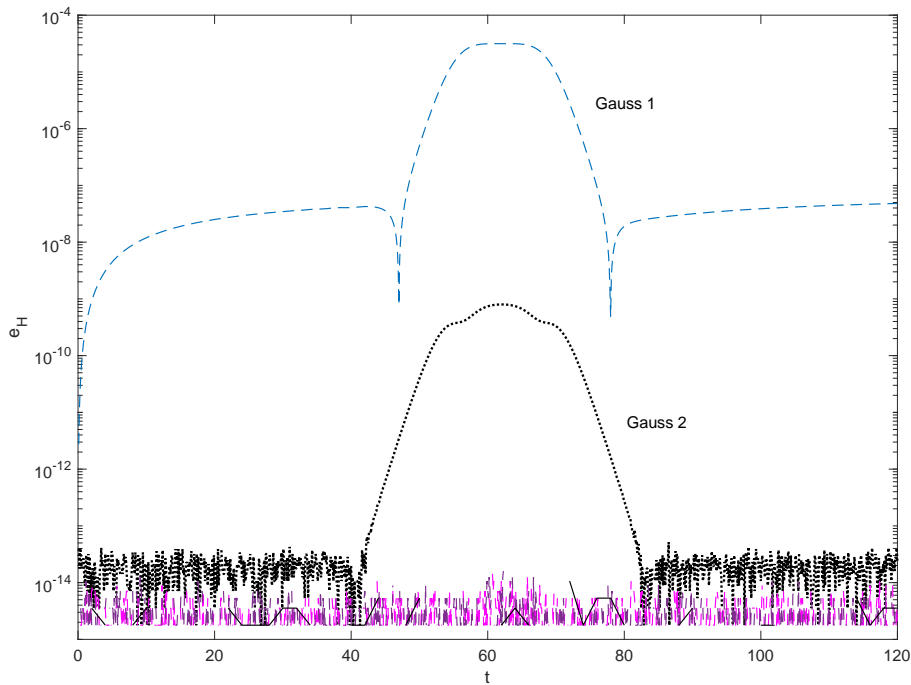}}
\caption{Hamiltonian error for the 1-stage Gauss method (dashed line), 2-stage Gauss method (dotted line), using a timestep $h=0.1$. The Hamiltonian errors close to round-off are those of HBVM(2,1) and HBVM(3,2), using a timestep $h=0.1$, and SHBVM, using timestep $h=2$.}
\label{eHfig}
\end{figure}

It is worth mentioning that, as is shown in Figure~\ref{eHfig}, for this problem the symplectic methods exhibit  a growth in the Hamiltonian error, when the two waves collide, unless the timestep is very small. Conversely, the energy conserving HBVMs and the SHBVM method  always provide a uniformly small Hamiltonian error.

\begin{table}[p]
\caption{Numerical results with timestep $h=120/n$ for the collision of solitary waves problem (\ref{uvform})-(\ref{cwuv0}).}
\label{cwtab}

\smallskip
\centerline{\begin{tabular}{|r|r|c|c|c|c|c|c|}
\hline
\multicolumn{8}{|c|}{Gauss 1}\\
\hline
$n$ & time & $e_u$ & rate & $e_H$ & rate & $e_M$ & ---\\ 
\hline
1200 &  7.6 & 7.62e-04 &  --- & 3.15e-05&  ---  & 1.60e-14 &   \\ 
2400 & 13.0 & 1.90e-04 &  2.0 & 7.87e-06&  2.0  & 1.59e-14 &   \\ 
3600 & 27.0 & 8.46e-05 &  2.0 & 3.50e-06&  2.0  & 1.60e-14 &  \\ 
4800 & 35.8 & 4.76e-05 &  2.0 & 1.97e-06&  2.0  & 1.60e-14 &  \\ 
6000 & 44.5 & 3.05e-05 &  2.0 & 1.26e-06&  2.0  & 1.59e-14 &   \\ 
\hline
\hline
\multicolumn{8}{|c|}{Gauss 2}\\
\hline
$n$ & time & $e_u$ & rate & $e_H$ & rate & $e_M$ & ---\\ 
\hline
1200 & 16.9 & 3.16e-08 &  --- & 7.98e-10&  ---  & 1.62e-14 &    \\ 
2400 & 30.3 & 1.97e-09 &  4.0 & 4.99e-11&  4.0  & 1.65e-14 &   \\ 
3600 & 40.7 & 3.90e-10 &  4.0 & 9.86e-12&  4.0  & 1.63e-14 &    \\ 
4800 & 49.4 & 1.23e-10 &  4.0 & 3.11e-12&  4.0  & 1.64e-14 &   \\ 
6000 & 61.8 & 5.05e-11 &  4.0 & 1.27e-12&  4.0  & 1.64e-14 &    \\ 
\hline
\hline
\multicolumn{8}{|c|}{HBVM(2,1)}\\
\hline
$n$ & time & $e_u$ & rate & $e_H$ & --- & $e_M$ & rate\\ 
\hline
1200 & 16.2 & 7.66e-04 &  ---  & 1.42e-14&   & 1.64e-14 &  --- \\ 
2400 & 28.9 & 1.91e-04 &  2.0 & 1.60e-14&   & 1.64e-14 &  ** \\ 
3600 & 38.8 & 8.51e-05 &  2.0 & 1.95e-14&   & 1.65e-14 & ** \\ 
4800 & 49.2 & 4.78e-05 &  2.0 & 1.78e-14&   & 1.64e-14 & ** \\ 
6000 & 64.1 & 3.06e-05 &  2.0 & 1.78e-14&   & 1.65e-14 & ** \\ 
\hline
\hline
\multicolumn{8}{|c|}{HBVM(3,2)}\\
\hline
$n$ & time & $e_u$ & rate & $e_H$ & --- & $e_M$ & rate\\ 
\hline
1200 & 19.4 & 3.13e-08 &  --- & 1.60e-14&     & 1.62e-14 &  --- \\ 
2400 & 34.3 & 1.96e-09 &  4.0 & 1.42e-14&     & 1.64e-14 & ** \\ 
3600 & 45.2 & 3.87e-10 &  4.0 & 1.95e-14&    & 1.64e-14 & ** \\ 
4800 & 60.9 & 1.22e-10 &  4.0 & 1.78e-14&     & 1.63e-14 &  ** \\ 
6000 & 67.8 & 5.02e-11 &  4.0 & 2.31e-14&    & 1.63e-14 & ** \\ 
\hline
\hline
\multicolumn{8}{|c|}{SHBVM $(k=18,s=12)$}\\
\hline
$n$ & time & $e_u$ & --- & $e_H$ & --- & $e_M$ & ---\\ 
\hline
60 & 11.4 & 7.86e-14 & & 1.07e-14 & & 1.59e-14 &\\
\hline
\end{tabular}}
\end{table}

\section{Conclusions}\label{end}
In this paper we have studied the efficient numerical solution of the ``good'' Boussinesq equation with periodic boundary conditions. The equation has, at first, been cast into Hamiltonian form, then using a spectrally accurate Fourier space discretization. Time integration has then been carried out by considering the energy conserving HBVM$(\lceil \frac{3}2s\rceil,s)$ methods. In particular, when $s$ is suitably large, such methods can be regarded as spectral methods in time (SHBVMs). A very efficient implementation of such methods, relying on their so-called {\em blended} implementation, has been then considered, providing a very efficient numerical method for solving the ``good'' Boussinesq equation, with spectral accuracy both in space and time. A few numerical tests duly confirm this conclusion, showing that SHBVMs provide, for the problem at hand, a geometric integrator able to preserve all the invariants of the problem, as well as to provide a negligible solution error. These results further confirm the effectiveness of SHBVMs for solving Hamiltonian PDEs \cite{BIMR2018}.

\section{Acknowledgements} 
This paper emerged during a visit of the third author in Firenze, which has been supported by NSFC (Grant No. 11571128) and from the Universit\`a di Firenze (ex 60\% project and ``Progetto di internazionalizzazione di Ateneo'' c/o DIMAI). 

The first author wishes also to thanks the chat ``{\em B\dots\/come Ben trovati}'' for the support during the nights spent in writing the paper. A particular thank is to Nino Losito, for the music selections.

\end{document}